\documentclass[11pt]{article}
\usepackage{amscd}
\usepackage{amsfonts}
\usepackage{amsmath}
\usepackage{amssymb}
\usepackage{amsthm}
\usepackage{bbm}
\usepackage{CJK}
\usepackage{fancyhdr}
\usepackage{graphicx}
\usepackage{hyperref}
\usepackage{indentfirst}
\usepackage{latexsym,bm}
\usepackage{mathrsfs}
\usepackage{xypic}

\usepackage[top=1in,bottom=1in,left=1.25in,right=1.25in]{geometry}
\def\<{\langle}
\def\>{\rangle}

\def\l{\ltimes}
\def\lr{\longrightarrow}

\def\o{\otimes}

\def\si{\sigma}

\date{}
\begin{document}
\renewcommand{\baselinestretch}{1.2}
\renewcommand{\arraystretch}{1.0}
\title{\bf The Drinfel'd Double versus the Heisenberg
 Double for Hom-Hopf Algebras }
 \date{}
\author {{\bf Daowei Lu, \quad  Shuanhong Wang \footnote {Corresponding author:  Shuanhong Wang, shuanhwang@seu.edu.cn, 
 Department of Mathematics, Southeast University, Nanjing, Jiangsu 210096, P. R. of China.}}\\
{\small Department of Mathematics, Southeast University}\\
{\small Nanjing, Jiangsu 210096, P. R. of China}}
 \maketitle
\begin{center}
\begin{minipage}{12.cm}

\noindent{\bf Abstract.} Let $(A,\alpha)$ be a finite-dimensional Hom-Hopf algebra. In this paper we mainly
 construct the Drinfel'd double $D(A)=(A^{op}\bowtie A^{\ast},\alpha\otimes(\alpha^{-1})^{\ast})$ in the setting of Hom-Hopf algebras
 by two ways, one of which
 generalizes  Majid's bicrossproduct for Hopf algebras (see \cite{M2})
 and another one is to introduce the notion of dual pairs of of Hom-Hopf algebras.
 Then we study the relation between the Drinfel'd double
 $D(A)$ and Heisenberg double  $H(A)=A\# A^{*}$, generalizing the main result in \cite{Lu}.
 Especially, the examples given in the paper are not obtained from the usual Hopf algebras.
 \\

\\

 \end{minipage}
 \end{center}
 \normalsize\vskip1cm

\section*{0. INTRODUCTION}

As well-known, in physics, the Yang-Baxter equation (or star-triangle relation) is a consistency equation
 which was first introduced in the field of statistical mechanics, and
 it takes its name from independent work of C. N. Yang in  1968, and R. J. Baxter in 1971.
 It depends on the idea that in some scattering situations, particles may preserve their
  momentum while changing their quantum internal states.
   It states that a matrix $R$, acting on two out of three objects,
   satisfies
   $$
   (R\o id)(id\o R)(R\o id)=(id\o R)(R\o id)(id\o R).
   $$
In one dimensional quantum systems, $R$ is the scattering matrix and if it satisfies the Yang-Baxter equation
 then the system is integrable. The Yang-Baxter equation also shows up when discussing knot theory
 and the braid groups where $R$ corresponds to swapping two strands.
 Since one can swap three strands two different ways, the Yang-Baxter equation enforces that both paths are the same.
\\

Braided monoidal  categories  give rise to solutions to the quantum Yang-Baxter equation (QYBE).
 The classical braided monoidal categories
 come from the representations of quasitriangular Hopf algebras.
  Quasitriangular Hopf algebras have been very widely studied and
 one of the most important examples of a quasitriangular Hopf algebra is the Drinfel'd double (see \cite{Dr}).
 As such, they are interesting to different research communities in mathematical physics (see \cite{M2}  for example).
  Although the Drinfel'd double for finite-dimensional Hopf algebras provide examples of such solutions to the QYBE,
  these are rather trivial.
  The wish to obtain more interesting solutions to the QYBE provides a strong
  motivation to find new examples of Drinfel'd double in the setting of Hom-Hopf algebras,
   whose representations are a braided monoidal category.
\\

The aim of this article is to construct new examples of braided monoidal categories in the setting of Hom-Hopf algebras.
 This is achieved by generalizing an existing the Drinfel'd double (see \cite{Dr}).
 We will generalize this construction to the so-called Hom-Hopf algebras,
  which are Hopf algebras in the Hom-category of a monoidal category (see \cite{CG}).
   We find a suitable generalization of the notions of a Majid's bicrossproduct and dual pairs
 for this setting and obtain a Drinfel'd double  of a Hom-Hopf algebra whose representation is a braided
 monoidal category. Furthermore,  we  establish  the relation between Drinfel'd double
 and Heisenberg double in the setting of Hom-Hopf algebras, generalizing the results in \cite{Lu}.
 Especially, the examples obtained in the paper are not obtained from the usual Hopf algebras.
 \\

The article is organized as follows.
\\

In Section 1, we will recall the definitions and results of  Hom-Hopf algebras, such as  Hom-algebas,  Hom-coalgebras,  Hom-modules,  Hom-comodules and the Hom-smash products.
\\

Let $(A,\alpha _A)$ and $(H,\alpha_H)$ be finite-dimensional Hom-Hopf algebras.
In Section 2, we will introduce the notion of a bicrossproduct $(A\l H,\alpha_{A}\o\alpha_{H})$, and give the conditions for $(A\l H,\alpha_{A}\o\alpha_{H})$ to form a Hom-Hopf algebra, generalizing the Majid's bicrossproduct defined in \cite{M1} (see  Theorem 2.6).
And examples of bicrossproduct  Hom-Hopf algebras are constructed. Then we will show how to get the Drinfel'd double by our new theory.
Furthermore solutions of the QHYBE are obtained.
\\

In Section 3, we will construct the Drinfel'd double associated to a dual pair of Hom-Hopf algebras (see Theorem 3.4),
 which is a generalization of the Drinfel'd double in the finite dimensional case.
 In section 4, we will establish the relation between Drinfel'd double and Heisenberg double (see Theorem 4.5).
 Examples are given to illustrate our theory.

\section*{1. PRELIMINARIES}
\def\theequation{1. \arabic{equation}}
\setcounter{equation} {0} \hskip\parindent

Throughout this article, all the vector spaces, tensor product and homomorphisms are over a fixed field $k$.  For a coalgebra $C$, we will use the Heyneman-Sweedler's notation $\Delta(c)=  c_{1}\otimes c_{2},$
for any $c\in C$ (summation omitted).
\\

In this section, we will recall the definitions in \cite{MP} on the  Hom-Hopf algebras,  Hom-modules and  Hom-comodules.
\\

 A unital Hom-associative algebra is a triple $(A,\mu,\alpha)$ where $\alpha:A\lr A$ and $\mu:A\o A\lr A$ are linear maps, with notation $\mu(a\o b)=ab$ such that for any $a,b,c\in A$,
 \begin{eqnarray*}
&&\alpha(ab)=\alpha(a)\alpha(b),\ \alpha(1_{A})=1_{A},\\
&&1_{A}a=\alpha(a)=a1_{A},\ \alpha(a)(bc)=(ab)\alpha(c).
\end{eqnarray*}
A linear map $f:(A,\mu_{A},\alpha_{A})\lr (B,\mu_{B},\alpha_{B})$ is called a morphism of Hom-associative algebra if $\alpha_{B}\circ f=f\circ\alpha_{A}$, $f(1_{A})=1_{B}$ and $f\circ\mu_{A}=\mu_{B}\circ(f\o f).$

 A counital Hom-coassociative coalgebra is a triple $(C,\Delta,\varepsilon,\alpha)$ where $\alpha:C\lr C$, $\varepsilon:C\lr k$, and $\Delta:C\lr C\o C$ are linear maps such that
 \begin{eqnarray*}
&&\varepsilon\circ\alpha=\varepsilon,\ (\alpha\o\alpha)\circ\Delta=\Delta\circ\alpha,\\
&&(\varepsilon\o id)\circ\Delta=\alpha=(id\o\varepsilon)\circ\Delta,\\
&&(\Delta\o\alpha)\circ\Delta=(\alpha\o\Delta)\circ\Delta.
\end{eqnarray*}

A linear map $f:(C,\Delta_{C},\alpha_{C})\lr (D,\Delta_{D},\alpha_{D})$ is called a morphism of Hom-coassociative coalgebra if $\alpha_{D}\circ f=f\circ\alpha_{C}$, $\varepsilon_{D}\circ f=\varepsilon_{C}$ and $\Delta_{D}\circ f=(f\o f)\circ\Delta_{C}.$

In what follows, we will always assume all Hom-algebras are unital and Hom-coalgebras are counital.

A Hom-bialgebra is a quadruple $(H,\mu,\Delta,\alpha)$, where $(H,\mu,\alpha)$ is a Hom-associative algebra and $(H,\Delta,\alpha)$ is a Hom-coassociative coalgebra such that $\Delta$ and $\varepsilon$ are morphisms of Hom-associative algebra.

A Hom-Hopf algebra $(H,\mu,\Delta,\alpha)$ is a Hom-bialgebra $H$ with a linear map $S:H\lr H$(called antipode) such that
\begin{eqnarray*}
&&S\circ\alpha=\alpha\circ S,\\
&&S(h_{1})h_{2}=  h_{1}S(h_{2})=\varepsilon(h)1,
\end{eqnarray*}
for any $h\in H$. For $S$ we have the following properties:
 \begin{eqnarray*}
&&S(h)_{1}\o S(h)_{2}=  S(h_{2})\o S(h_{1}),\\
&&S(gh)=S(h)S(g),\ \varepsilon\circ S=\varepsilon.
 \end{eqnarray*}
For any Hopf algebra $H$ and any Hopf algebra endomorphism $\alpha$ of $H$, there exists a Hom-Hopf algebra $H_{\alpha}=(H,\alpha\circ\mu,1_{H},\Delta\circ\alpha,\varepsilon,S,\alpha)$.

Let $(A,\alpha_{A})$ be a Hom-associative algebra, $M$ a linear space and $\alpha_{M}:M\lr M$ a linear map. A left $A$-module structure on $(M,\alpha_{M})$ consists of a linear map $A\o M\lr M$, $a\o m\mapsto a\cdot m$, such that
 \begin{eqnarray*}
&&1_{A}\cdot m=\alpha_{M}(m),\\
&&\alpha_{M}(a\cdot m)=\alpha_{A}(a)\cdot\alpha_{M}(m),\\
&&\alpha_{A}(a)\cdot(b\cdot m)=(ab)\cdot\alpha_{M}(m),
\end{eqnarray*}
for any $a,b\in A$ and $m\in M.$

Similarly we can define the right $(A,\alpha)$-modules. Let $(M,\mu)$ and $(N,\nu)$ be two left $(A,\alpha)$-modules, then a linear map $f:M\lr N$ is a called left $A$-module map if $f(am)=af(m)$ for any $a\in A$, $m\in M$ and $f\circ\mu=\nu\circ f$.

Let $(C,\alpha_{C})$ be a Hom-coassociative coalgebra, $M$ a linear space and $\alpha_{M}:M\lr M$ a linear map. A right $C$-comodule structure on $(M,\alpha_{M})$ consists a linear map $\rho:M\lr M\o C$ such that
 \begin{eqnarray*}
&&(id\o\varepsilon_{C})\circ\rho=\alpha_{M},\\
&&(\alpha_{M}\o\alpha_{C})\circ\rho=\rho\circ\alpha_{M},\\
&&(\rho\o\alpha_{C})\circ\rho=(\alpha_{M}\o\Delta)\circ\rho.
 \end{eqnarray*}
Let $(M,\mu)$ and $(N,\nu)$ be two right $(C,\gamma)$-comodules, then a linear map $g:M\lr N$ is a called right $C$-comodule map if $g\circ \mu=\nu\circ g$ and $\rho_{N}\circ g=(g\otimes id)\circ\rho_{M}$.

Let $(H,\mu_{H},\Delta_{H},\alpha_{H})$ be a Hom-bialgebra. A Hom-associative algebra $(A,\mu_{A},\alpha_{A})$ is called a left $H$-module Hom-algebra if $(A,\alpha_{A})$ is a left $H$-module, with the action $H\o A\lr A,\ h\o a\mapsto h\cdot a$, such that
$$\alpha^{2}_{H}(h)\cdot (ab)=  (h_{1}\cdot a)(h_{2}\cdot b),$$
$$h\cdot 1_{A}=\varepsilon(h)1_{A},$$
for any $h\in H$ and $a,b\in A$.

When $A$ is a left $H$-module Hom-algebra, in \cite{MP} the Hom-smash product $A\#H$ is defined as follows:
$$(a\#h)(b\#k)=  a(\alpha^{-2}_{H}(h_{1})\cdot\alpha^{-1}_{A}(b))\#\alpha^{-1}_{H}(h_{2})k,$$
for any $a,b\in A$ and $h,k\in H$.

Recall from
\cite{Y2} that
a Hom-bialgebra $(H,\alpha_{H})$ is quasitriangular if there exists an element $R\in H\otimes H$ satisfying:
\begin{enumerate}
\item
$\Delta^{op}(x)R=R\Delta(x)$ for any $x\in H,$
\item
$(\Delta\otimes\alpha_{H})R=R^{13}R^{23},$
\item
$(\alpha_{H}\otimes\Delta)R=R^{13}R^{12}.$
\end{enumerate}

\section*{2. MAJID'S BICROSSPRODUCTS FOR HOM-HOPF ALGEBRAS}
\def\theequation{2. \arabic{equation}}
\setcounter{equation} {0} \hskip\parindent

In this section, we will construct the bicrossproduct on Hom-bialgebras, which generalize
the Majid's bicrossproducts in the usual Hopf algebras. Then we can get Drinfel'd double for
Hom-Hopf algebras in the setting of the Majid's bicrossproducts.
\\

In \cite{LS} L. Liu and B. Shen defined the smash coproduct on monoidal Hom-Hopf algebra. Now we need to construct the smash coproduct on Hom-Hopf algebra.
\\

{\bf Definition 2.1.}
Let $(H,\alpha_{H})$ be a Hom-bialgebra. A Hom-coassociative coalgebra $(C,\alpha_{C})$ is called a right {\sl $H$-comodule Hom-coalgebra} if $(C,\alpha_{C})$ a right $H$-comodule, with the comodule structure map $\rho:C\lr C\o H$, $c\mapsto   c_{(0)}\o c_{(1)}$, such that the following conditions hold
$$ \varepsilon_{C}(c_{(0)})c_{(1)}=\varepsilon_{C}(c)1_{H},$$
$$  c_{(0)1}\o c_{(0)2}\o\alpha^{2}_{H}(c_{(1)})=  c_{1(0)}\o c_{2(0)}\o c_{1(1)}c_{2(1)},$$
for any $c\in C.$
\\

{\bf Example 2.2.} (1) Let $H$ be a bialgebra and $C$ a right $H$-comodule coalgebra in the usual sense. The coaction is denoted by $C\lr C\o H$, $c\mapsto  c_{[0]}\o c_{[1]}$. Let $\alpha_{H}$ be a bialgebra endomorphism of $H$, and $\alpha_{C}$ a coalgebra endomorphism of $C$ such that
$$ \alpha_{C}(c)_{[0]}\o \alpha_{C}(c)_{[1]}= \alpha_{C}(c_{[0]})\o\alpha_{H}(c_{[1]}),$$
for any $c\in C.$ Consider the Hom-bialgebra $H_{\alpha_{H}}=(H,\alpha_{H}\circ\mu_{H},\Delta_{H}\circ\alpha_{H},\alpha_{H})$ and the Hom-coassociative coalgebra $C_{\alpha_{C}}=(C,\Delta_{C}\circ\alpha_{C},\alpha_{C})$. Then $C_{\alpha_{C}}$ is a right $H_{\alpha_{H}}$-comodule Hom-coalgebra with the coaction $C_{\alpha_{C}}\lr C_{\alpha_{C}}\o H_{\alpha_{H}}$, $c\mapsto   c_{(0)}\o c_{(1)}=  \alpha_{C}(c_{[0]})\o \alpha_{H}(c_{[1]}).$

(2) Let $(H,\alpha)$ be any Hom-Hopf algebra. Then $(H^{op},\alpha)$ is a right $(H,\alpha)$-comodule Hom-coalgebra with the right  coaction
$$\rho(h)=  \alpha^{-1}(h_{12})\o S(\alpha^{-2}(h_{11}))\alpha^{-1}(h_{2}),$$
for any $h\in H$.  The fact that $(H^{op},\alpha)$ is a right $(H,\alpha)$-comodule is obvious, and for $h\in H$,
$$\begin{aligned}
&  h_{(0)1}\o h_{(0)2}\o \alpha^{2}(h_{1})\\
=&  \alpha^{-1}(h_{121})\o\alpha^{-1}(h_{122})\o S(h_{11})\alpha(h_{2})\\
=&  h_{12}\o \alpha^{-1}(h_{212})\o[S(\alpha^{-1}(h_{11}))(\alpha^{-4}(h_{2111})S(\alpha^{-4}(h_{2112})))]h_{22}\\
=&  \alpha^{-1}(h_{112})\o \alpha^{-1}(h_{212})\o[S(\alpha^{-2}(h_{111}))(\alpha^{-2}(h_{12})S(\alpha^{-3}(h_{211})))]h_{22}\\
=&  \alpha^{-1}(h_{112})\o \alpha^{-1}(h_{212})\o[S(\alpha^{-2}(h_{111}))\alpha^{-1}(h_{12})][S(\alpha^{-2}(h_{211}))\alpha^{-1}(h_{22})]\\
=&  h_{1(0)}\o h_{2(0)}\o h_{1(1)}h_{2(1)}.
\end{aligned}$$
Then $(H^{op},\alpha)$ is a right $(H,\alpha)$-comodule Hom-coalgebra.

(3)  Let $A^1_x=span\{1,x\}$ over a fixed field $k$ with $chark\neq2$. Then define $\beta$ as a $k$-linear automorphism of $A^1_x$ by
$$\beta(1_{A^1_x})=1_{A^1_x},\ \ \beta(x)=-x.$$
Define the multiplication on $A^1_x$ by
$$1_{A^1_x}1_{A^1_x}=1_{A^1_x},\ 1_{A^1_x}x=-x,\ x^{2}=0,$$
then it is not hard to check that  $(A^1_x,\beta)$ is a Hom-associative algebra.

 For $A^1_x$,  define the coalgebra structure and antipode by
$$\Delta(1_{A^1_x})=1_{A^1_x}\otimes 1_{A^1_x},\ \Delta(x)=(-x)\otimes 1+1\otimes(-x),$$
$$\varepsilon(1_{A^1_x})=1,\ \varepsilon(x)=0,\ \ S_{A^1_x}(1_{A^1_x})=1_{A^1_x},\ S_{A^1_x}(x)=-x,$$

then $(A^1_x,\beta)$ is a Hom-Hopf algebra.

Let $H=span\{1,g|g^{2}=1\}$ be the group algebra. Obviously $(H,id)$ is a Hom-Hopf algebra.
Define the left action of $H$ on $A^1_x$ $\cdot:H\otimes A^1_x\lr A^1_x$ by
$$1_{H}\cdot 1_{A^1_x}=1_{A^1_x},\ 1_{H}\cdot x=-x,\ g\cdot 1_{A^1_x}=1_{A^1_x},\ g\cdot x=x.$$
It is easy to check that $(A^1_x,\beta)$ is a left $(H,id)$-module algebra.

Define the right coaction of $A^1_x$ on $H$ $\rho:H\lr H\otimes A^1_x$ by
$$\rho(1_{H})=1_{H}\otimes 1_{A^1_x},\ \rho(g)=g\otimes 1_{A^1_x},$$
then $(H,id)$ is a right $(A^1_x,\beta)$-comodule coalgebra.
\\

{\bf Definition 2.3.}
Let $(C,\Delta_{C},\alpha_{C})$ and $(D,\Delta_{D},\alpha_{D})$ be Hom-coassociative coalgebras. A linear map $\Phi:C\o D\lr D\o C$ is called a {\sl Hom-cotwisting map} between $C$ and $D$ if the following conditions hold:
\begin{eqnarray}
&&  (\Delta_{D}\o \alpha_{C})\circ \Phi=(id_{D}\o\Phi)\circ(\Phi\o id_{D})\circ(\alpha_{C}\o\Delta_{D}),\\
&& (\alpha_{D}\o\Delta_{C})\circ \Phi=(\Phi\o id_{C})\circ(id_{C}\o\Phi)\circ(\Delta_{C}\o\alpha_{D}), \\
&& (\alpha_{D}\o\alpha_{C})\circ\Phi=\Phi\circ(\alpha_{C}\o\alpha_{D}),\\
&& (\varepsilon_{C}\o id)\circ \Phi=id\o\varepsilon_{C},\ (id\o\varepsilon_{D})\circ \Phi=\varepsilon_{D}\o id.
\end{eqnarray}

If we denote the element $\Phi(c\o d)$ in $D\o C$ by $\Phi(c\o d)=  d^{\Phi}\o c^{\Phi}=  d^{\phi}\o c^{\phi}$ for any $c\in C,d\in D$, the above identities can be rewritten as
\[\left\{
\begin{array}{ccccccc}
  &(d^{\Phi})_{1}\o(d^{\Phi})_{2}\o\alpha_{C}(c^{\Phi})=  (d_{1})^{\Phi}\o(d_{2})^{\phi}\o\alpha_{C}(c)^{\Phi\phi},\\
  &\alpha_{D}(d^{\Phi})\o(c^{\Phi})_{1}\o(c^{\Phi})_{2}= \alpha_{D}(d)^{\Phi\phi}\o(c_{1})^{\phi}\o(c_{2})^{\Phi},\\
  &\alpha_{D}(d^{\Phi})\o\alpha_{C}(c^{\Phi})= \alpha_{D}(d)^{\Phi}\o\alpha_{C}(c)^{\Phi},\\
 &\varepsilon_{C}(c^{\Phi})d^{\Phi}=\varepsilon_{C}(c)d,\  \varepsilon_{D}(d^{\Phi})c^{\Phi}=\varepsilon_{D}(d)c.
\end{array} \right.
\]
\\

{\bf Proposition 2.4.}
Let $(C,\Delta_{C},\alpha_{C})$ and $(D,\Delta_{D},\alpha_{D})$ be Hom-coassociative coalgebras and  $\Phi:C\o D\lr D\o C$ a Hom-cotwisting map. Define $\Delta:C\o D\lr (C\o D)\o(C\o D)$ by
$$\Delta(c\o d)=  c_{1}\o d_{1}^{\Phi}\o c_{2}^{\Phi}\o d_{2},$$
and
$\varepsilon:C\o D\lr k$ by
$$\varepsilon(c\o d)=\varepsilon_{C}(c)\varepsilon_{D}(d),$$
for any $c\in C,d\in D$. Then $(C\o D, \Delta, \alpha_{C}\o\alpha_{D})$ is a Hom-coassociative coalgebra.

{\bf Proof.} The counit is straightforward to verify.
For any $c\in C,d\in D$, firstly
$$\begin{aligned}
(\Delta\circ(\alpha_{C}\o\alpha_{D}))(c\o d)&=  \Delta(\alpha_{C}(c)\o\alpha_{D}(d))\\
                                          &= \alpha_{C}(c_{1})\o\alpha_{D}(d_{1})^{\Phi}\o\alpha_{C}(c_{2})^{\Phi}\o\alpha_{D}(d_{2})\\
                                          &= \alpha_{C}(c_{1})\o\alpha_{D}(d_{1}^{\Phi})\o\alpha_{C}(c_{2}^{\Phi})\o\alpha_{D}(d_{2})\\
                                          &=((\alpha_{C}\o\alpha_{D})\o(\alpha_{C}\o\alpha_{D}))\circ\Delta(c\o d).
\end{aligned}$$

Then
\begin{eqnarray*}
&&((\alpha_{C}\o\alpha_{D})\o\Delta)\circ\Delta(c\o d)\\
&=&\alpha_{C}(c_{1})\o\alpha_{D}(d_{1}^{\Phi})\o\Delta(c_{2}^{\Phi}\o d_{2})\\
&=&\alpha_{C}(c_{1})\o\alpha_{D}(d_{1}^{\Phi})\o(c_{2}^{\Phi})_{1}\o(d_{21})^{\phi}\o((c_{2}^{\Phi})_{2})^{\phi}\o d_{22}\\
&\stackrel{(2.1)}{=}& \alpha_{C}(c_{1})\o\alpha_{D}(d_{1})^{\Phi\phi'}\o(c_{21})^{\phi'}\o(d_{21})^{\phi}\o(c_{22})^{\Phi\phi}\o d_{22}\\
&=&c_{11}\o (d_{11})^{\Phi\phi'}\o(c_{12})^{\phi'}\o(d_{12})^{\phi}\o\alpha_{C}(c_{2})^{\Phi\phi}\o \alpha_{D}(d_{2})\\
&\stackrel{(2.2)}{=}&  c_{11}\o ((d_{1}^{\Phi})_{1})^{\phi}\o(c_{12})^{\phi}\o(d_{1}^{\Phi})_{2}\o\alpha_{C}(c_{2}^{\Phi})\o \alpha_{D}(d_{2})\\
&=&(\Delta\o(\alpha_{C}\o\alpha_{D}))\circ\Delta(c\o d).
\end{eqnarray*}

By definition, this finishes the proof. \hfill $\blacksquare$
\\

{\bf Example  2.5.} (1) Let $(H,\alpha_{H})$ be a Hom-bialgebra and $(C,\alpha_{C})$ a right $H$-comodule Hom-coalgbra with the coaction $C\lr C\o H,\ c\mapsto c_{(0)}\o c_{(1)}$. Define the linear map $\Phi_{1}: H\o C\lr C\o H$ by
$$\Phi_{1}(h\o c)=\alpha^{-1}_{C}(c_{(0)})\o\alpha^{-1}_{H}(h)\alpha^{-2}_{C}(c_{(1)}),
$$
for any $h\o H$ and $c\in C$. Then $\Phi_{1}$ is a Hom-cotwisting map between $H$ and $C$.

(2) Let $(H,\alpha)$ be a Hom-bialgebra and $\rho$ the coaction of $H$ on $H^{op}$ in Example 2.2 (2), we have a Hom-cotwisting map $\Phi_{2}:H\o H^{op}\lr H^{op}\o H$ given by
$$\Phi_{2}(k\o h)=\alpha^{-2}(h_{12})\o\alpha^{-1}(k)[S(\alpha^{-4}(h_{11}))\alpha^{-3}(h_{2})],
$$
for any $h,\ k\in H$.

(3) With the notations in Example 2.2 (3), we have a Hom-cotwisting map $\Phi_{3}:A^{1}_{x}\o H\lr H\o A^{1}_{x}$ given by
\begin{eqnarray*}
&&\Phi_{3}(1_{A^1_x}\o1_{H})=1_{H}\o1_{A^1_x},\quad \Phi_{3}(x\o1_{H})=1_{A^1_x}\o-x,\\
&&\Phi_{3}(1_{A^1_x}\o g)=g\o 1_{A^1_x},\quad\quad \Phi_{3}(x\o g)=g\o-x.
\end{eqnarray*}
\\

{\bf Theorem 2.6.} Let $(A,S_{A},\alpha_{A})$ and $(H,S_{H},\alpha_{H})$ be Hom-Hopf algebras. Assume that $A$ is a left $H$-module Hom-algebra with the action $H\o A\lr A,\ h\o a\mapsto h\cdot a$ and $H$ is a right $A$-comodule Hom-coalgebra with the coaction $H\lr H\o A,\ h\mapsto  h_{(0)}\o h_{(1)}$. If the following conditions are satisfied:

(1) For any $h\in H,b\in A$,
$$\Delta(h\cdot b)= \alpha^{-1}_{H}(h_{1(0)})\cdot b_{1}\o \alpha^{-1}_{A}(h_{1(1)})(\alpha^{-1}_{H}(h_{2})\cdot \alpha^{-1}_{A}(b_{2})),$$
$$\varepsilon_{A}(h\cdot b)=\varepsilon_{A}(b)\varepsilon_{H}(h),$$

(2) For any $h,g\in H$,
$$ (hg)_{(0)}\o(hg)_{(1)}= \alpha^{-1}_{H}(h_{1(0)})g_{(0)}\o\alpha^{-1}_{A}(h_{1(1)})(\alpha^{-1}_{H}(h_{2})\cdot\alpha^{-1}_{A}(g_{(1)})),$$

(3) For any $h\in H,b\in A$,
$$  h_{2(0)}\o (h_{1}\cdot b)h_{2(1)}=  h_{1(0)}\o h_{1(1)}(h_{2}\cdot b),$$
then $(A\o H, \alpha_{A}\o\alpha_{H})$ is a Hom-Hopf algebra under the Hom-smash product and Hom-smash coproduct, that is,
$$(a\o h)(b\o g)=  a(\alpha^{-2}_{H}(h_{1})\cdot\alpha^{-1}_{A}(b))\o\alpha^{-1}_{H}(h_{2})g,$$
and
$$\Delta(a\o h)=  a_{1}\o\alpha^{-1}_{H}(h_{1(0)})\o\alpha^{-1}_{A}(a_{2})\alpha^{-2}_{A}(h_{1(1)})\o h_{2},$$
with the antipode $S:A\o H\lr A\o H$ given by
$$S(a\o h)= (1\o S_{H}(\alpha^{-2}_{H}(h_{(0)})))(S_{A}(\alpha^{-2}_{A}(a)\alpha^{-3}_{A}(h_{(1)}))\o1),$$
for any $h,g\in H$, $a,b\in A$.

We will call this Hom-Hopf algebra Hom-bicrossproduct, i.e, {\sl Majid's bicrossproduct} for
 Hom-Hopf algebras, and denote it by $(A\ltimes H, \alpha_{A}\o\alpha_{H})$.

{\bf Proof.} Firstly we need to prove $\Delta$ and $\varepsilon$ is a morphism of Hom-associative algebra. Indeed for any $a,b\in A$ and $g,h\in H$,
$$\begin{aligned}
&\Delta((a\o h)(b\o g))\\
=&\Delta(  a(\alpha^{-2}_{H}(h_{1})\cdot\alpha^{-1}_{A}(b))\o\alpha^{-1}_{H}(h_{2})g)\\
=&  a_{1}(\alpha^{-2}_{H}(h_{1})\cdot\alpha^{-1}_{A}(b)_{1})\o\alpha^{-1}_{H}((\alpha^{-1}_{H}(h_{21})g_{1})_{(0)})\\
&\o[\alpha^{-1}_{A}(a_{2})\alpha^{-1}_{A}((\alpha^{-2}_{H}(h_{1})\cdot\alpha^{-1}_{A}(b))_{2})]\alpha^{-2}_{A}
((\alpha^{-1}_{H}(h_{21})g_{1})_{(1)})\o\alpha^{-1}_{H}(h_{22})g_{2}\\
=&  a_{1}(\alpha^{-3}_{H}(h_{11(0)})\cdot\alpha^{-1}_{A}(b_{1}))\o\alpha^{-3}_{H}(h_{211(0)})\alpha^{-1}_{H}(g_{1(0)})\\
&\o\{\alpha^{-1}_{A}(a_{2})[\alpha^{-4}_{A}(h_{11(1)})(\alpha^{-4}_{H}(h_{12})\cdot\alpha^{-3}_{A}(b_{2}))]\}\\
&\{\alpha^{-4}_{A}(h_{211(1)})(\alpha^{-4}_{H}(h_{212})\cdot\alpha^{-3}_{A}(g_{1(1)}))\}\o\alpha^{-1}_{H}(h_{22})g_{2}\\
=&  a_{1}(\alpha^{-3}_{H}(h_{11(0)})\cdot\alpha^{-1}_{A}(b_{1}))\o\alpha^{-3}_{H}(h_{211(0)})\alpha^{-1}_{H}(g_{1(0)})\\
&\o\{\alpha^{-1}_{A}(a_{2})\alpha^{-3}_{A}(h_{11(1)})\}\{[(\alpha^{-3}_{H}(h_{12})\cdot\alpha^{-2}_{A}(b_{2}))]\\
&[\alpha^{-5}_{A}(h_{211(1)})(\alpha^{-5}_{H}(h_{212})\cdot\alpha^{-4}_{A}(g_{1(1)}))]\}\o\alpha^{-1}_{H}(h_{22})g_{2}\\
\end{aligned}$$

$$\begin{aligned}
=&  a_{1}(\alpha^{-3}_{H}(h_{11(0)})\cdot\alpha^{-1}_{A}(b_{1}))\o\alpha^{-3}_{H}(h_{211(0)})\alpha^{-1}_{H}(g_{1(0)})\o\{\alpha^{-1}_{A}(a_{2})\alpha^{-3}_{A}(h_{11(1)})\}\\
&\{[(\alpha^{-4}_{H}(h_{12})\cdot\alpha^{-3}_{A}(b_{2}))\alpha^{-5}_{A}(h_{211(1)})](\alpha^{-4}_{H}(h_{212})\cdot\alpha^{-3}_{A}(g_{1(1)}))\}\o\alpha^{-1}_{H}(h_{22})g_{2}\\
=&  a_{1}(\alpha^{-2}_{H}(h_{1(0)})\cdot\alpha^{-1}_{A}(b_{1}))\o\alpha^{-4}_{H}(h_{2112(0)})\alpha^{-1}_{H}(g_{1(0)})\o\{\alpha^{-1}_{A}(a_{2})\alpha^{-2}_{A}(h_{1(1)})\}\\
&\{[(\alpha^{-6}_{H}(h_{2111})\cdot\alpha^{-3}_{A}(b_{2}))\alpha^{-6}_{A}(h_{2112(1)})](\alpha^{-4}_{H}(h_{212})\cdot\alpha^{-3}_{A}(g_{1(1)}))\}\o\alpha^{-1}_{H}(h_{22})g_{2}\\
=&  a_{1}(\alpha^{-2}_{H}(h_{1(0)})\cdot\alpha^{-1}_{A}(b_{1}))\o\alpha^{-4}_{H}(h_{2111(0)})\alpha^{-1}_{H}(g_{1(0)})\o\{\alpha^{-1}_{A}(a_{2})\alpha^{-2}_{A}(h_{1(1)})\}\\
&\{[\alpha^{-6}_{A}(h_{2111(1)})(\alpha^{-6}_{H}(h_{2112})\cdot\alpha^{-3}_{A}(b_{2}))](\alpha^{-4}_{H}(h_{212})\cdot\alpha^{-3}_{A}(g_{1(1)}))\}\o\alpha^{-1}_{H}(h_{22})g_{2}\\
=&  a_{1}(\alpha^{-2}_{H}(h_{1(0)})\cdot\alpha^{-1}_{A}(b_{1}))\o\alpha^{-4}_{H}(h_{2111(0)})\alpha^{-1}_{H}(g_{1(0)})\o\{\alpha^{-1}_{A}(a_{2})\alpha^{-2}_{A}(h_{1(1)})\}\\
&\{\alpha^{-5}_{A}(h_{2111(1)})[(\alpha^{-6}_{H}(h_{2112})\cdot\alpha^{-3}_{A}(b_{2}))(\alpha^{-5}_{H}(h_{212})\cdot\alpha^{-4}_{A}(g_{1(1)}))]\}\o\alpha^{-1}_{H}(h_{22})g_{2}\\
=&  a_{1}(\alpha^{-2}_{H}(h_{1(0)})\cdot\alpha^{-1}_{A}(b_{1}))\o\alpha^{-3}_{H}(h_{211(0)})\alpha^{-1}_{H}(g_{1(0)})\o\{\alpha^{-1}_{A}(a_{2})\alpha^{-2}_{A}(h_{1(1)})\}\\
&\{\alpha^{-4}_{A}(h_{211(1)})[\alpha^{-4}_{H}(h_{212})\cdot\alpha^{-3}_{A}(b_{2})\alpha^{-4}_{A}(g_{1(1)})]\}\o\alpha^{-1}_{H}(h_{22})g_{2}\\
=&  a_{1}(\alpha^{-2}_{H}(h_{1(0)})\cdot\alpha^{-1}_{A}(b_{1}))\o\alpha^{-3}_{H}(h_{211(0)})\alpha^{-1}_{H}(g_{1(0)})\\
&\o\{\alpha^{-1}_{A}(a_{2})[\alpha^{-3}_{A}(h_{1(1)})\alpha^{-5}_{A}(h_{211(1)})]\}\{\alpha^{-3}_{H}(h_{212})\cdot\alpha^{-2}_{A}(b_{2})\alpha^{-3}_{A}(g_{1(1)})\}\o\alpha^{-1}_{H}(h_{22})g_{2}\\
=&  a_{1}(\alpha^{-4}_{H}(h_{111(0)})\cdot\alpha^{-1}_{A}(b_{1}))\o\alpha^{-3}_{H}(h_{112(0)})\alpha^{-1}_{H}(g_{1(0)})\\
&\o\{\alpha^{-1}_{A}(a_{2})\alpha^{-5}_{A}(h_{111(1)}h_{112(1)})\}\{\alpha^{-2}_{H}(h_{12})\cdot\alpha^{-2}_{A}(b_{2})\alpha^{-3}_{A}(g_{1(1)})\}\o h_{2}g_{2}\\
=&  a_{1}(\alpha^{-4}_{H}(h_{11(0)1})\cdot\alpha^{-1}_{A}(b_{1}))\o\alpha^{-3}_{H}(h_{11(0)2})\alpha^{-1}_{H}(g_{1(0)})\o\{\alpha^{-1}_{A}(a_{2})\alpha^{-3}_{A}(h_{11(1)})\}\\
&\{\alpha^{-2}_{H}(h_{12})\cdot\alpha^{-2}_{A}(b_{2})\alpha^{-3}_{A}(g_{1(1)})\}\o h_{2}g_{2}\\
=&  a_{1}(\alpha^{-3}_{H}(h_{1(0)1})\cdot\alpha^{-1}_{A}(b_{1}))\o\alpha^{-2}_{H}(h_{1(0)2})\alpha^{-1}_{H}(g_{1(0)})\o\{\alpha^{-1}_{A}(a_{2})\alpha^{-2}_{A}(h_{1(1)})\}\\
&\{\alpha^{-2}_{H}(h_{21})\cdot\alpha^{-2}_{A}(b_{2})\alpha^{-3}_{A}(g_{1(1)})\}\o \alpha^{-1}_{H}(h_{22})g_{2}\\
=&\Delta(a\o h)\Delta(b\o g).
\end{aligned}$$

It is straightforward to check that $\varepsilon$ is a morphism of Hom-associate algebra.

For the antipode $S$ we have
$$\begin{aligned}
& (a\o h)_{1}S((a\o h)_{2})\\
=& (a_{1}\o\alpha^{-1}_{H}(h_{1(0)}))[(1\o S_{H}(\alpha^{-2}_{H}(h_{2(0)}))\\
&(S_{A}((\alpha^{-3}_{A}(a_{2})\alpha^{-4}_{A}(h_{1(1)}))\alpha^{-3}_{A}(h_{2(1)}))\o 1)]\\
=& [a_{1}\o\alpha^{-2}_{H}(h_{1(0)})S_{H}(\alpha^{-2}_{H}(h_{2(0)}))]\\
&[S_{A}(\alpha^{-1}_{A}(a_{2})(\alpha^{-3}_{A}(h_{1(1)})\alpha^{-3}_{A}(h_{2(1)})))\o 1]\\
=& [a_{1}\o\alpha^{-2}_{H}(h_{(0)1})S_{H}(\alpha^{-2}_{H}(h_{(0)2}))][S_{A}(\alpha^{-1}_{A}(a_{2})(\alpha^{-1}_{A}(h_{(1)})))\o 1]\\
=& [a_{1}\o\varepsilon(h_{(0)})1][S_{A}(\alpha^{-1}_{A}(a_{2})(\alpha^{-1}_{A}(h_{(1)})))\o 1]\\
=&\varepsilon_{A}(a)\varepsilon_{H}(h)1\o1.
\end{aligned}$$

Similarly we have $  S(a\o h)_{1}(a\o h)_{2}=\varepsilon_{A}(a)\varepsilon_{H}(h)1\o1.$ Thus $A\o H$ is a Hom-Hopf algebra.

This finishes the proof. \hfill $\blacksquare$
\\

{\bf Example 2.7.} We now consider  $A^1_x$ given in Example 2.2 (3).
By a simple computation the conditions in Theorem 2.6 are satisfied.
 Hence we have a bicrossproduct Hom-Hopf algebra $(A\#H,\beta\#id)$ with the
 following Hom-Hopf algebra structure:

 The product is given by:
 $$\begin{array}{|c|c|c|c|c|}
\hline                                & 1_{A}\#1_{H} & 1_{A}\#g & x\#1_{H} & x\# g \\
\hline 1_{A}\#1_{H}                  & 1_{A}\#1_{H} & 1_{A}\#g & -x\#1_{H} & -x\# g \\
       1_{A}\#g                      & 1_{A}\#g & 1_{A}\#1_{H} & x\#g & x\#1_{H} \\
       x\#1_{H}                      & -x\#1_{H} & -x\#g &  0 & 0      \\
       x\# g                         &  -x\# g  & -x\# 1 & 0 & 0         \\
 \hline
 \end{array}$$

The coproduct is defined as:
$$\Delta(1_{A}\#1_{H})=1_{A}\#1_{H}\o1_{A}\#1_{H},\ \Delta(1_{A}\#g)=1_{A}\#g\o1_{A}\#g,$$
$$\Delta(x\#1_{H})=-x\#1_{H}\o1_{A}\#1_{H}-1_{A}\#1_{H}\o x\#1_{H},$$
$$\Delta(x\o g)=1_{A}\#g\o x\# g-x\# g\o1_{A}\#g,$$

and the antipode is given by:
$$S(1_{A}\#1_{H})=1_{A}\#1_{H},\ \ S(1_{A}\#g )=1_{A}\#g,$$
$$S(x\#1_{H})=-x\#1_{H},\ \ S(x\# g)=x\# g.$$
\\

{\bf Example 2.8.}
Let $A$ and $H$ be two Hopf algebras, and $A\#H$ the bicrossproduct of $A$ and $H$ with $A$ a left $H$-module algebra and $H$ a right $A$-comodule coalgebra. Assume that $\alpha_{A}$ and $\alpha_{H}$ are the automorphisms of Hopf algebras of $A$ and $H$ respectively, and satisfy
\begin{eqnarray*}
&&\alpha_{A}(h\cdot a)=\alpha_{H}(h)\cdot\alpha_{A}(a),\\
&&\alpha_{H}(h)_{(0)}\o\alpha_{H}(h)_{(1)}=\alpha_{H}(h_{(0)})\o\alpha_{A}(h_{(1)}).
\end{eqnarray*}
Then $A_{\alpha_{A}}$ is a left $H_{\alpha_{H}}$-module algebra under the action $h\triangleright a=\alpha_{H}(h)\cdot\alpha_{A}(a)$ and $H_{\alpha_{H}}$ is a right $A_{\alpha_{A}}$-comodule coalgebra under the coaction $h_{[0]}\o h_{[1]}=\alpha_{H}(h_{(0)})\o\alpha_{A}(h_{(1)})$. Hence we have a Hom-type bicrossproduct $A_{\alpha_{A}}\#H_{\alpha_{H}}$.
\\

Next we will construct a class of Hom-Hopf algebras.
\\

{\bf Corollary 2.9.} For any Hom-Hopf algebra $(H,S,\alpha)$, there exists a bicrossproduct structure on the space $(H\o H^{op},\alpha\o\alpha)$.

{\bf Proof.} From the Example 2.5 and Proposition 2.4, $(H^{op},\alpha)$ is a right $(H,\alpha)$-comodule Hom-coalgebra. And $(H,\alpha)$ is a left $(H^{op},\alpha)$-module Hom-algebra under the action
$$h\cdot a= (S(\alpha^{-2}(h_{1}))\alpha^{-1}(a))\alpha^{-1}(h_{2}),$$
for any $a,h\in H$. The verification of module condition is straightforward and is left to the reader.

For any $h,a,b\in H$,
$$\begin{aligned}
&  (h_{1}\cdot a)(h_{2}\cdot b)\\
=& [(S(\alpha^{-2}(h_{11}))\alpha^{-1}(a))\alpha^{-1}(h_{12})][(S(\alpha^{-2}(h_{21}))\alpha^{-1}(b))\alpha^{-1}(h_{22})]\\
=& [(S(\alpha^{-2}(h_{11}))\alpha^{-1}(a))(\alpha^{-2}(h_{12})S(\alpha^{-2}(h_{21}))][b\alpha^{-1}(h_{22})]\\
=& [(S(\alpha^{-1}(h_{1}))\alpha^{-1}(a))(\alpha^{-3}(h_{211})S(\alpha^{-2}(h_{212}))][b\alpha^{-1}(h_{22})]\\
=&  (S(h_{1})a)(bh_{2})=  (S(h_{1})\alpha^{-1}(ab))\alpha(h_{2})\\
=&\alpha^{2}(h)\cdot (ab).
\end{aligned}$$
That is, $(H,\alpha)$ is a left $(H^{op},\alpha)$-module Hom-algebra.

Moreover for any $h,g,a\in H$,
$$\begin{aligned}
&  (h\cdot a)_{1}\o(h\cdot a)_{2}\\
=&  (S(\alpha^{-2}(h_{12}))\alpha^{-1}(a_{1}))\alpha^{-1}(h_{21})\o(S(\alpha^{-2}(h_{11}))\alpha^{-1}(a_{2}))\alpha^{-1}(h_{22})\\
=&  (S(\alpha^{-3}(h_{121}))\alpha^{-1}(a_{1}))\alpha^{-2}(h_{122})\o[S(\alpha^{-2}(h_{11}))(\alpha^{-4}(h_{211})S\alpha^{-4}(h_{212}))]\\
&(\alpha^{-1}(a_{2})\alpha^{-2}(h_{22}))\\
=&  (S(\alpha^{-4}(h_{1121}))\alpha^{-1}(a_{1}))\alpha^{-3}(h_{1122})\o[S(\alpha^{-3}(h_{111}))\alpha^{-2}(h_{12})]\\
&[(S\alpha^{-3}(h_{21})\alpha^{-2}(a_{2}))\alpha^{-2}(h_{22})]\\
=&  (S(\alpha^{-3}(h_{1(0)1}))\alpha^{-1}(a_{1}))\alpha^{-2}(h_{1(0)2})\o\alpha^{-1}(h_{1(1)})[(S\alpha^{-3}(h_{21})\alpha^{-2}(a_{2}))\alpha^{-2}(h_{22})]\\
=&  \alpha^{-1}(h_{1(0)}\cdot a_{1}\o \alpha^{-1}(h_{1(1)})(\alpha^{-1}(h_{2})\cdot\alpha^{-1}(b_{2})),
\end{aligned}$$

and
$$\begin{aligned}
& (gh)_{0}\o(gh)_{1}= \alpha^{-1}(g_{12}h_{12})\o S(\alpha^{-2}(g_{11}h_{11}))\alpha^{-1}(g_{2}h_{2})\\
=& \alpha^{-1}(g_{12})\alpha^{-1}(h_{12})\o [S(\alpha^{-2}(h_{11}))S(\alpha^{-2}(g_{11}))][\alpha^{-1}(g_{2})\alpha^{-1}(h_{2})]\\
=& \alpha^{-1}(g_{12})\alpha^{-1}(h_{12})\o S(\alpha^{-1}(h_{11}))[(S(\alpha^{-3}(g_{11}))\alpha^{-2}(g_{2}))\alpha^{-1}(h_{2})]\\
=& \alpha^{-1}(g_{12})\alpha^{-2}(h_{112})\o [S(\alpha^{-3}(h_{111}))\alpha^{-2}(h_{12})]\\
&[(S(\alpha^{-3}(h_{21}))(S(\alpha^{-4}(g_{11}))\alpha^{-3}(g_{2})))\alpha^{-2}(h_{22})]\\
=&  g_{(0)}\alpha^{-1}(h_{1(0)})\o\alpha^{-1}(h_{1(1)})(\alpha^{-1}(h_{2})\cdot\alpha^{-1}(g_{(1)})),
\end{aligned}$$

Finally
$$\begin{aligned}
&  h_{2(0)}\o (h_{1}\cdot a)h_{2(1)}\\
=& \alpha^{-1}(h_{212})\o[(S(\alpha^{-2}(h_{11}))\alpha^{-1}(a))\alpha^{-1}(h_{12})][S(\alpha^{-2}(h_{211}))\alpha^{-1}(h_{22})]\\
=& \alpha^{-1}(h_{212})\o[(S(\alpha^{-2}(h_{11}))\alpha^{-1}(a))(\alpha^{-2}(h_{12})S(\alpha^{-3}(h_{211})))]h_{22}\\
=& \alpha^{-1}(h_{122})\o[(S(\alpha^{-2}(h_{11}))\alpha^{-1}(a))(\alpha^{-4}(h_{1211})S(\alpha^{-4}(h_{1212})))]\alpha(h_{2})\\
=&  h_{12}\o(S(\alpha^{-1}(h_{11}))a)\alpha(h_{2})\\
=&  h_{12}\o S(h_{11})(ah_{2})\\
=&  h_{12}\o [S(\alpha^{-1}(h_{11}))(\alpha^{-3}(h_{211})S(\alpha^{-3}(h_{212}))](a\alpha^{-1}(h_{22}))\\
=&  h_{1(0)}\o h_{1(1)}(h_{2}\cdot a).
\end{aligned}$$
This finishes the proof. \hfill $\blacksquare$
\\

Therefore by Theorem 2.6, we have a bicrossproduct structure on $H\o H^{op}$ with multiplication and comultiplication as follows:
$$(a\ltimes h)(b\ltimes k)=  a[(S(\alpha^{-4}(h_{11}))\alpha^{-2}(b))\alpha^{-3}(h_{12})]\ltimes k\alpha^{-1}(h_{2}),$$
$$\Delta(a\ltimes h)=  a_{1}\ltimes\alpha^{-2}(h_{112})\o\alpha^{-1}(a_{2})(S(\alpha^{-4}(h_{111}))\alpha^{-3}(h_{12}))\ltimes h_{2},$$
for any $a,b,h,k\in H$.
\\

{\bf Definition 2.10.}
Let $(H,\alpha_{H})$ be a Hom-bialgebra and $(C,\alpha_{C})$ a Hom-coassociative coalgebra. Then $C$ is called a left {\sl $H$-module Hom-coalgebra} if $C$ is  a left $H$-module with the action $H\o C\lr C,\ h\o c\mapsto h\cdot c$, such that
$$\Delta_{C}(h\cdot c)=  h_{1}\cdot c_{1}\o h_{2}\cdot c_{2},\ \ \varepsilon_{C}(h\cdot c)=\varepsilon_{H}(h)\varepsilon_{C}(c),$$
for any $h\in H,\ c\in C.$
\\

{\bf Definition 2.11.}
Let $(A,\alpha_{A})$ and $(H,\alpha_{H})$ be two Hom-bialgebras. $(A,H)$ is called a {\sl matched pair} if there exist linear maps
$$\triangleleft:H\o A\lr H,\ \triangleright:H\o A\lr A,$$
turning $A$ into a left $H$-module Hom-coalgebra and turning $H$ into a right $A$-module Hom-coalgebra such that the following conditions are satisfied:
\begin{eqnarray}
&&(hg)\triangleleft a= (h\triangleleft(\alpha^{-2}_{H}(g_{1})\triangleright\alpha^{-3}_{A}(a_{1}))(\alpha^{-1}_{H}(g_{2})\triangleleft\alpha^{-2}_{A}(a_{2}))\\
&&h\triangleright(ab)= (\alpha^{-2}_{H}(h_{1})\triangleright\alpha^{-1}_{A}(a_{1}))((\alpha^{-3}_{H}(h_{2})\triangleleft\alpha^{-2}_{A}(a_{2}))\triangleright b),\\
&&h_{1}\triangleleft a_{1}\o h_{2}\triangleright a_{2}=  h_{2}\triangleleft a_{2}\o h_{1}\triangleright a_{1}.
\end{eqnarray}
\\

{\bf Proposition 2.12.}
Let $(A,S_{A},\alpha_{A})$ and $(H,S_{H},\alpha_{H})$ be two Hom-Hopf algebras, and $(A,H)$ a matched pair. There exists a unique Hom-Hopf algebra structure on the vector space $A\o H$ with the multiplication, comultiplication and  the antipode given by
\begin{eqnarray*}
&&(a\o h)(b\o g)=  a(\alpha^{-2}_{H}(h_{1})\triangleright\alpha^{-2}_{A}(b_{1}))\o(\alpha^{-2}_{H}(h_{2})\triangleleft\alpha^{-2}_{A}(b_{2}))g,\\
&&\Delta(a\o h)=  a_{1}\o h_{1}\o a_{2}\o h_{2},\\
&&S(a\o h)=(1_{A}\o S_{H}\alpha^{-1}_{H}(h))(S_{A}\alpha^{-1}_{A}(a)\o 1_{H}),
\end{eqnarray*}
for any $a,b\in A$ and $g,h\in H$.

Equipped with this Hopf algebra structure, $A\o H$ is called double crossed product of $A$ and $H$ denoted by $A\bowtie H$.

{\bf Proof.} Define the linear map $T:H\o A\lr A\o H$ by
$$T(h\o a)= \alpha^{-2}_{H}(h_{1})\triangleright\alpha^{-2}_{A}(b_{1})\o\alpha^{-2}_{H}(h_{2})\triangleleft\alpha^{-2}_{A}(b_{2}),$$
for any $a\in A$ and $h\in H$.

By Proposition 2.6 in \cite{MP}, in order to prove the multiplication is Hom-associative,  we need only to verify that $T$ is a Hom-twisting map between $H$ and $A$. Indeed, first easy to see that
$$T\circ(\alpha_{H}\o\alpha_{A})=(\alpha_{A}\o\alpha_{H})\circ T.$$

Then for any $a,b\in A$ and $h\in H$,
$$\begin{aligned}
& (ab)_{T}\o\alpha_{H}(h)_{T}\\
=& \alpha^{-1}_{H}(h_{1})\triangleright\alpha^{-2}_{A}(a_{1}b_{1})\o\alpha^{-1}_{H}(h_{2})\triangleleft\alpha^{-2}_{A}(a_{2}b_{2})\\
=& \alpha^{-1}_{H}(h_{1})\triangleright\alpha^{-2}_{A}(a_{1}b_{1})\o\alpha^{-1}_{H}(h_{2})\triangleleft\alpha^{-2}_{A}(a_{2}b_{2})\\
=& (\alpha^{-3}_{H}(h_{11})\triangleright\alpha^{-3}_{A}(a_{11}))((\alpha^{-4}_{H}(h_{12})\triangleleft\alpha^{-4}_{A}(a_{12}))\triangleright\alpha^{-2}(b_{1}))\o\alpha^{-1}_{H}(h_{2})\triangleleft\alpha^{-2}_{A}(a_{2}b_{2})\\
=& (\alpha^{-2}_{H}(h_{1})\triangleright\alpha^{-2}_{A}(a_{1}))((\alpha^{-4}_{H}(h_{21})\triangleleft\alpha^{-4}_{A}(a_{21}))\triangleright\alpha^{-2}(b_{1}))\\
&\o(\alpha^{-3}_{H}(h_{22})\triangleleft\alpha^{-3}_{A}(a_{22}))\triangleleft\alpha^{-1}_{A}(b_{2})\\
=&  a_{T}b_{t}\o\alpha_{H}(h_{Tt}),
\end{aligned}$$
where we have used the notation $T(h\o a)=a_{T}\o h_{T}=a_{t}\o h_{t}$.

Similarly for any $a\in A$ and $g,h\in H$, we have
$$ \alpha_{A}(a)_{T}\o(hg)_{T}= \alpha_{A}(a_{Tt})\o h_{t}g_{T}.$$ Now $A\o H$ is a Hom-associative algebra. Next for any $a,b\in A,h,g\in H$,

$$\begin{aligned}
&\Delta((a\o h)(b\o g))\\
=& \Delta(a(\alpha^{-2}_{H}(h_{1})\triangleright\alpha^{-2}_{A}(b_{1}))\o \alpha^{-2}_{H}(h_{2})\triangleleft\alpha^{-2}_{A}(b_{2}))g\\
=&  a_{1}(\alpha^{-2}_{H}(h_{11})\triangleright\alpha^{-2}_{A}(b_{11}))\o (\alpha^{-2}_{H}(h_{21})\triangleleft\alpha^{-2}_{A}(b_{21}))g_{1}\\
&\o a_{2}(\alpha^{-2}_{H}(h_{12})\triangleright\alpha^{-2}_{A}(b_{12}))\o (\alpha^{-2}_{H}(h_{22})\triangleleft\alpha^{-2}_{A}(b_{22}))g_{2}\\
=&  a_{1}(\alpha^{-2}_{H}(h_{11})\triangleright\alpha^{-2}_{A}(b_{11}))\o (\alpha^{-2}_{H}(h_{12})\triangleleft\alpha^{-2}_{A}(b_{12}))g_{1}\\
&\o a_{2}(\alpha^{-2}_{H}(h_{21})\triangleright\alpha^{-2}_{A}(b_{21}))\o (\alpha^{-2}_{H}(h_{22})\triangleleft\alpha^{-2}_{A}(b_{22}))g_{2}\\
=&\Delta(a\o h)\Delta(b\o g).
\end{aligned}$$
Therefore $A\o H$ is a Hom-bialgebra. It is straightforward to check that $S$ is the antipode.

This finishes the proof. \hfill $\blacksquare$
\\

Let $(H,\alpha)$ be a Hom-bialgebra.
In \cite{MP}, the authors have defined the Hom-associative algebra $(H^{*},(\alpha^{-1})^{*})$, the linear dual of $H$, where the multiplication is given by
$$(f\bullet g)(h)=  f(\alpha^{-2}(h_{1}))g(\alpha^{-2}(h_{2})),$$
for any $f,g\in H^{*}$, and $h\in H$.

Now we give the comultiplication on $H^{*}$ by
$$f(gh)=  f_{1}(\alpha^{2}(h))f_{2}(\alpha^{2}(g)),$$
for any $f\in H^{*}$ and $g,h\in H$. In other words, $\langle\Delta(f),h\o k\rangle=f(\alpha^{-2}(hk))$.
\\

{\bf Proposition 2.13.}
Let $(H,S,\alpha)$ be a Hom-Hopf algebra. With the multiplication and comultiplication defined on $H^{*}$ as above, $(H^{*},S^{*},(\alpha^{-1})^{*})$ is a Hom-Hopf algebra.

{\bf Proof}~~First we need to prove that $(H^{*},(\alpha^{-1})^{*})$ is a Hom-coassociative coalgebra. Indeed for any $g,h,k\in H$ and $f\in H^{*}$,
$$\begin{aligned}
f((gh)\alpha(k))&=  f_{1}(\alpha^{2}(gh))f_{2}(\alpha^{3}(k))\\
                &=  f_{11}(\alpha^{4}(g))f_{12}(\alpha^{4}(h))f_{2}(\alpha^{3}(k)),
\end{aligned}$$
and
$$\begin{aligned}
f(\alpha(g)(hk))&=  f_{1}(\alpha^{3}(g))f_{2}(\alpha^{2}(hk))\\
                &=  f_{1}(\alpha^{3}(g))f_{21}(\alpha^{4}(h))f_{22}(\alpha^{4}(k)),
\end{aligned}$$
therefore $  (\alpha^{-1})^{*}(f_{1})\o f_{21}\o f_{22}=  f_{11}\o f_{12}\o(\alpha^{-1})^{*}(f_{2}).$

Furthermore
$$\langle\Delta\circ(\alpha^{-1})^{*}(f),h\o k\rangle=(\alpha^{-1})^{*}(f)(\alpha^{-2}(hk))=f(\alpha^{-3}(hk)),$$
and
$$ (\alpha^{-1})^{*}(f_{1})(h)(\alpha^{-1})^{*}(f_{2})(k)=  f_{1}(\alpha^{-1}(h))f_{2}(\alpha^{-1}(k))=f(\alpha^{-3}(hk)).$$
That is, $\Delta\circ(\alpha^{-1})^{*}=((\alpha^{-1})^{*}\o(\alpha^{-1})^{*})\circ\Delta$. Then $(H^{*},(\alpha^{-1})^{*})$ is a Hom-coassociative coalgebra.

Finally for any $f,g\in H^{*}$ and $h,k\in H$,
$$\begin{aligned}
\langle\Delta(f\bullet g),h\o k\rangle&=(f\bullet g)(\alpha^{-2}(hk))\\
                                      &=  f(\alpha^{-4}(h_{1}k_{1}))g(\alpha^{-4}(h_{2}k_{2}))\\
                                      &=  f_{1}(\alpha^{-2}(h_{1}))f_{2}(\alpha^{-2}(k_{1}))g_{1}(\alpha^{-2}(h_{2}))g_{2}(\alpha^{-2}(k_{2}))\\
                                      &=  \langle f_{1}\bullet g_{1}\o f_{2}\bullet g_{2},h\o k\rangle.
\end{aligned}$$
Therefore $(H^{*},(\alpha^{-1})^{*})$ is a Hom-bialgebra. It is not hard to verify that $S^{*}$ is the antipode.

This finishes the proof. \hfill $\blacksquare$
\\

The following proposition can be found in \cite{MP}.
\\

{\bf Proposition 2.14.} Let $(C,\alpha_{C})$ be a Hom-coassociative coalgebra, and $(M,\alpha_{M})$ is a right $C$-comodule. Then $(M,\alpha_{M})$ is a left $(C^{*},(\alpha^{-1}_{C})^{*})$-module with the action $f\cdot m=  f(m_{(1)})m_{(0)}$.
\\

{\bf Proposition 2.15.}
Let $(A,\alpha_{A})$ and $(H,\alpha_{H})$ be two Hom-bialgebras, $(A\ltimes H,\alpha_{A}\o\alpha_{H})$ the bicrossproduct of $A$ and $H$. Define the left action $\triangleright:A^{*}\o H\lr H$ of $A^{*}$ on $H$ by
$$f\triangleright h=  f(h_{(1)})h_{(0)},$$
and the right action $\triangleleft:A^{*}\o H\lr A^{*}$ of $H$ on $A^{*}$ by
$$\langle f\triangleleft h,a\rangle=\langle f,h\cdot\alpha^{-2}(a)\rangle,$$
for any $f\in A^{*}$, $h\in H$ and $a\in A$. Then $(A^{*},(\alpha^{-1}_{A})^{*})$ and $(H,\alpha_{H})$ is a matched pair.

{\bf Proof.} First, we need to verify that $H$ is a left $A^{*}$-module Hom-coalgebra.  By Proposition 2.14, $H$ is a left $A^{*}$-module. Then for any $h\in H$ and $f\in A^{*}$,
$$\begin{aligned}
\Delta(f\triangleright h)&=  f(h_{(1)})h_{(0)1}\o h_{(0)2}\\
                         &=  f(\alpha^{-2}(h_{1(1)}h_{2(1)}))h_{1(0)}\o h_{2(0)}\\
                         &=  f_{1}(h_{1(1)})h_{1(0)}\o f_{2}(h_{2(1)})h_{2(0)}\\
                         &=  f_{1}\triangleright h_{1}\o f_{2}\triangleright h_{2}.
\end{aligned}$$
So $H$ is a left $A^{*}$-module Hom-coalgebra.

Now for any $g,h\in H$, $f\in A^{*}$ and $a,b\in A$
$$\begin{aligned}
\langle(\alpha^{-1}_{A})^{*}(f)\triangleleft(hg),a\rangle&=\langle f,\alpha^{-1}_{H}(hg)\cdot\alpha^{-3}_{A}(a)\rangle\\
                                                         &=\langle f\triangleleft h,\alpha_{H}(g)\cdot\alpha^{-2}_{A}(a)\rangle\\
                                                         &=\langle (f\triangleleft h)\triangleleft\alpha_{H}(g),a\rangle.
\end{aligned}$$
and obviously
$$(\alpha^{-1}_{A})^{*}(f\triangleleft h)=(\alpha^{-1}_{A})^{*}(f)\triangleleft\alpha_{H}(h).$$

Then
$$\begin{aligned}
\langle\Delta(f\triangleleft h),a\o b\rangle&=\langle f\triangleleft h,\alpha^{-2}_{A}(ab)\rangle=\langle f,h\cdot\alpha^{-4}_{A}(ab)\rangle\\
                                            &= \langle f,(\alpha^{-2}_{H}(h_{1})\cdot\alpha^{-4}_{A}(a))(\alpha^{-2}_{H}(h_{2})\cdot\alpha^{-4}_{A}(b))\rangle\\
                                            &= \langle f_{1}, h_{1}\cdot\alpha^{-2}_{A}(a)\rangle\langle f_{2}, h_{2}\cdot\alpha^{-2}_{A}(b)\rangle\\
                                            &= \langle f_{1}\triangleleft h_{1},a\rangle\langle f_{2}\triangleleft h_{2},b\rangle.
\end{aligned}$$
Therefore $A^{*}$ is a right $H$-module Hom-coalgebra.

Next we will verify the compatible conditions in Definition 2.11. For any $f,g\in A^{*}$, $a\in A$ and $h,k\in H$,
$$\begin{aligned}
&\langle(f\bullet g)\triangleleft h,\ a\rangle=\langle f\bullet g,\ h\cdot\alpha^{-2}_{A}(a)\rangle\\
=&\langle f,\ \alpha^{-2}_{A}((h\cdot\alpha^{-2}_{A}(a))_{1})\rangle\langle  g,\ \alpha^{-2}_{A}((h\cdot\alpha^{-2}_{A}(a))_{2})\rangle\\
=&\langle f,\ \alpha^{-3}_{H}(h_{1(0)})\cdot\alpha^{-4}_{A}(a_{1})\rangle\langle  g,\ \alpha^{-3}_{H}(h_{1(1)})(\alpha^{-3}_{H}(h_{2})\cdot\alpha^{-5}_{A}(a_{2}))\rangle\\
=&\langle f,\ \alpha^{-3}_{H}(h_{1(0)})\cdot\alpha^{-4}_{A}(a_{1})\rangle\langle  g_{1},\ \alpha^{-1}_{H}(h_{1(1)})\rangle\langle g_{2},\ \alpha^{-1}_{H}(h_{2})\cdot\alpha^{-3}_{A}(a_{2})\rangle\\
=&\langle f,\ ((\alpha^{*}_{A})^{2}(g_{1})\triangleright\alpha^{-3}_{H}(h_{1}))\cdot\alpha^{-4}_{A}(a_{1})\rangle\langle \alpha_{A}^{*}(g_{2}),\ \alpha^{-2}_{H}(h_{2})\cdot\alpha^{-4}_{A}(a_{2})\rangle\\
=&\langle f\triangleleft((\alpha^{*}_{A})^{2}(g_{1})\triangleright\alpha^{-3}_{H}(h_{1})),\ \alpha^{-2}_{A}(a_{1})\rangle\langle \alpha_{A}^{*}(g_{2})\triangleleft\alpha^{-2}_{H}(h_{2}),\ \alpha^{-2}_{A}(a_{2})\rangle\\
=&\langle[f\triangleleft((\alpha^{*}_{A})^{2}(g_{1})\triangleright\alpha^{-3}_{H}(h_{1}))]\bullet[\alpha_{A}^{*}(g_{2})\triangleleft\alpha^{-2}_{H}(h_{2})], a\rangle,
\end{aligned}$$
thus we have the condition (2.5). And by
$$\begin{aligned}
f\triangleright(hk)&=  f((hk)_{(1)})(hk)_{(0)}\\
&= \langle f,\ \alpha^{-1}_{A}(h_{1(1)})(\alpha^{-1}_{H}(h_{2})\cdot\alpha^{-1}_{A}(k_{(1)}))\rangle\alpha^{-1}_{H}(h_{1(0)})k_{(0)}\\
&= \langle f_{1},\ \alpha_{A}(h_{1(1)})\rangle\langle f_{2},\ \alpha_{H}(h_{2})\cdot\alpha_{A}(k_{(1)})\rangle\alpha^{-1}_{H}(h_{1(0)})k_{(0)}\\
&= \langle (\alpha^{*}_{A})^{2}(f_{1}),\ \alpha^{-1}_{A}(h_{1(1)})\rangle\langle (\alpha^{*}_{A})^{3}(f_{2})\triangleleft \alpha^{-2}_{H}(h_{2}),\ k_{(1)}\rangle\alpha^{-1}_{H}(h_{1(0)})k_{(0)}\\
&= [(\alpha^{*}_{A})^{2}(f_{1})\triangleright\alpha^{-1}_{H}(h_{1})][((\alpha^{*}_{A})^{3}(f_{2})\triangleleft \alpha^{-2}_{H}(h_{2}))\triangleright k],
\end{aligned}$$
we get the condition (2.6). As for the condition (2.7),
$$\begin{aligned}
 \langle f_{1}\triangleleft h_{1},\ a\rangle\langle g,\ f_{2}\triangleright h_{2}\rangle&= \langle f_{1},\ h_{1}\cdot \alpha^{-2}_{A}(a)\rangle\langle g,\ h_{2(0)}\rangle\langle f_{2},\ h_{2(1)}\rangle\\
&= \langle f,\ \alpha^{-2}_{A}((h_{1}\cdot \alpha^{-2}_{A}(a))h_{2(1)})\rangle\langle g,\ h_{2(0)}\rangle\\
&= \langle f,\ \alpha^{-2}_{A}(h_{1(1)}(h_{2}\cdot \alpha^{-2}_{A}(a)))\rangle\langle g,\ h_{1(0)}\rangle\\
&= \langle f_{1},\ h_{1(1)}\rangle\langle f_{2},\ h_{2}\cdot \alpha^{-2}_{A}(a)\rangle\langle g,\ h_{1(0)}\rangle\\
&=  \langle f_{2}\triangleleft h_{2},\ a\rangle\langle g,\ f_{1}\triangleright h_{1}\rangle.
\end{aligned}$$
Therefore $(A^{*},(\alpha^{-1}_{A})^{*})$ and $(H,\alpha_{H})$ is a matched pair.

 This finishes the proof. \hfill $\blacksquare$
\\

By this result, we have the double crossed product $(H\bowtie A^{*},\alpha_{H}\o (\alpha^{-1}_{A})^{*})$ with the multiplication and comultiplication given by
$$(h\o f)(k\o g)=  h\alpha^{-2}_{H}(k_{1(0)})\o\langle f,\ \alpha^{-2}_{A}(k_{1(1)})(\alpha^{-2}_{H}(k_{2})\cdot\alpha^{-2}_{A}(?))\rangle g,$$
$$\Delta(h\o f)=  h_{1}\o f_{1}\o h_{2}\o f_{2},$$
for any $h,k\in H$ and $f,g\in A^{*}.$
\\

{\bf Corollary 2.16.} Let $(H,\alpha)$ be a Hom-Hopf algebra. We have the Drinfeld double $(H^{op}\bowtie H^{*},\alpha\o(\alpha^{-1})^{*})$ with the tensor comultiplication and multiplication given by
$$(h\o f)(k\o g)= \alpha^{-2}(k_{21})h\o[\alpha^{-3}(k_{22})\rightharpoonup((\alpha^{*})^{2}(f)\leftharpoonup S\alpha^{-3}(k_{1}))]g,$$
for any $h,k\in H$ and $f,g\in H^{*}$, where $\langle f\leftharpoonup h,\ k\rangle=\langle f, \ h\alpha^{-2}(k)\rangle$ and $\langle h\rightharpoonup f ,\ k\rangle=\langle f, \ \alpha^{-2}(k)h\rangle$.
\\

{\bf Example 2.17.} (1) Let $G$ be a finite group and $\phi$ an automorphism of $G$. Then $(kG,\phi)$ with the following structure is a Hom-Hopf algebra:
$$g\cdot h=\phi(gh),\ \Delta(g)=\phi(g)\otimes \phi(g),\ \varepsilon(g)=1,\ S(g)=g^{-1}.$$

Let $\{e_{g}\}_{g\in G}$ be the dual basis of the basis of $kG$. Then we have the Hom-Hopf algebra $k^{G}$, dual of $kG$, with the multiplication
$$e_{g}e_{h}=\delta_{g,h}\ e_{\phi(g)},$$
and the comultiplication, counit and antipode
$$\Delta^{*}(e_{g})= \sum_{uv=\phi(g)}e_{u}\o e_{v},\ \varepsilon^{*}(e_{g})=\delta_{g,1},\ S^{*}(e_{g})=e_{g^{-1}},$$
for any $g,h\in G.$

By Corollary 2.16, the multiplication in $D(kG)$ is given by
$$(g\o e_{h})(p\o e_{q})=\sum\phi(pg)\o \delta_{\phi(p)h\phi(p^{-1}),q}\ e_{\phi(q)}.$$

(2) In particular let $G=(g)$ be a cycle group of order $n>0$. Define $\phi\in Aut(G)$ by
$$\phi(g^{i})=g^{-i},\quad 0\leq i< n.
$$
Obviously $\phi\in Aut(G).$ Then we have a Hom-Hopf algebra $(kG,\phi)$ by
$$g^{i}\cdot g^{j}=g^{n-(i+j)},\ \Delta(g^{i})=g^{n-i}\o g^{n-i},\ \varepsilon(g^{i})=1,\ S(g^{i})=g^{n-i},
$$
for any $0\leq i,j< n$. Let $\{e_{i}\}_{0\leq i< n}$ be the dual basis of $kG$ such that $e^{i}(g^{j})=\delta_{i,j}e_{i}$. Then
$$\Delta^{*}(e_{i})= \sum_{j+k=n-i}e_{j}\o e_{k},\ \varepsilon^{*}(e_{i})=\delta_{i,1},\ S^{*}(e_{i})=e_{n-i}.$$
In the Drinfeld double $D(kG)$, we have the product
$$(g^{i}\o e_{j})(g^{m}\o e_{k})=\sum g^{n-(i+m)}\o \delta_{j,k}\ e_{n-k},$$
for any $0\leq i,j,m,k< n$.
\\

{\bf Example 2.18.} Let $H$ be the Hom-algebra generated by the elements $1_{H},\ g$ and $x$ satisfying the following relations:
$$1_{H}1_{H}=1_{H},\ 1_{H}g=g1_{H}=g,\ 1_{H}x=x1_{H}=-x,$$
$$g^{2}=1_{H},\ x^{2}=0,\ gx=-xg.$$
The automorphism $\alpha:H\lr H$ is defined by
$$\alpha(1_{H})=1_{H},\ \alpha(g)=g,\ \alpha(x)=-x,\ \alpha(gx)=-gx.$$
Then $(H,\alpha)$ is a Hom-associative algebra, and $\alpha^{2}=id$.

Define
$$\Delta(1_{H})=1_{H}\otimes 1_{H},\ \Delta(g)=g\otimes g,$$
$$\Delta(x)=(-x)\otimes g+1\otimes (-x),$$
$$\varepsilon(1_{H})=1,\ \varepsilon(g)=1,\ \varepsilon(x)=0,$$
$$S(1_{H})=1_{H},\ S(g)=g,\ S(x)=-gx.$$
Then $(H,\alpha)$ is a monoidal Hom-Hopf algebra.

Furthermore, let $R=\frac{1}{2}(1\otimes 1+1\otimes g+g\otimes 1-g\otimes g)$. Then
  it is not hard to check that $(H,\alpha,R)$ is quasitriangular.
\\

{\bf Proposition 2.19.}  The Drinfeld double $(H^{op}\bowtie H^{\ast},\alpha\otimes(\alpha^{-1})^{\ast})$ has the quasitriangular structure
$$R=  1\bowtie(\alpha^{-1})^{*}(e^{i})\o S^{-1}(e_{i})\bowtie\varepsilon,$$
where $\{e_{i}\}$ and $\{e^{*}_{i}\}$ are a base of $H$ and its dual base in $H^{\ast}$ respectively.

{\bf Proof.} For any $h,k,k'\in H$ and $f,g,g'\in H^{\ast}$, on one hand,
$$\begin{aligned}
&\langle\Delta^{op}(h\o f)R, g\o k\o g'\o k'\rangle\\
=& \langle (h_{2}\o f_{2})(1\o(\alpha^{-1})^{*}(e^{i})),g\o k\rangle\langle (h_{1}\o f_{1})(S^{-1}(e_{i})\o\varepsilon),g'\o k'\rangle\\
=& \langle g,\alpha(h_{2})\rangle\langle f_{2},\alpha^{-2}(k_{1})\rangle\langle g',\ S^{-1}(\alpha^{-5}(k_{2 21}))h_{1}\rangle\\
&\langle  f_{1},\alpha^{-5}(k_{222}))(\alpha^{-3}(k')S^{-1}(\alpha^{-5}(k_{21}))\rangle\\
=& \langle g,\alpha(h_{2})\rangle\langle g',\ S^{-1}(\alpha^{-4}(k_{21}))h_{1}\rangle\\
&\langle  f_{1},\alpha^{-4}(k_{22}))(\alpha^{-3}(k')S^{-1}(\alpha^{-5}(k_{12}))\rangle\langle f_{2},\alpha^{-3}(k_{11})\rangle\\
=& \langle g,\alpha(h_{2})\rangle\langle g',\ S^{-1}(\alpha^{-4}(k_{21}))h_{1}\rangle\\
&\langle  f,[\alpha^{-6}(k_{22}))(\alpha^{-5}(k')S^{-1}(\alpha^{-7}(k_{12}))]\alpha^{-5}(k_{11})\rangle\\
=& \langle g,\alpha(h_{2})\rangle\langle g',\ S^{-1}(\alpha^{-4}(k_{21}))h_{1}\rangle\\
&\langle  f,[\alpha^{-6}(k_{22})\alpha^{-4}(k')][S^{-1}(\alpha^{-6}(k_{12}))\alpha^{-6}(k_{11})]\rangle\\
=& \langle g,\alpha(h_{2})\rangle\langle g',\ S^{-1}(\alpha^{-3}(k_{1}))h_{1}\rangle\langle  f,[\alpha^{-4}(k_{2})\alpha^{-3}(k')]\rangle.
\end{aligned}$$

On the other hand,
$$\begin{aligned}
&\langle R\Delta(h\o f), g\o k\o g'\o k'\rangle\\
=& \langle (1\o(\alpha^{-1})^{*}(e^{i}))(h_{1}\o f_{1}),g\o k\rangle\langle (S^{-1}(\alpha(e_{i}))\o\varepsilon)(h_{2}\o f_{2}),g'\o k'\rangle\\
=& \langle g,\alpha^{-1}(h_{112})\rangle\langle e^{i},S\alpha^{-3}(h_{111})(\alpha^{-5}(k_{1})\alpha^{-3}(h_{12}))\rangle\langle f_{1},\alpha^{-2}(k_{2})\rangle\\
&\langle g',h_{2}S^{-1}(\alpha(e_{i}))\rangle\langle f_{2},\alpha^{-1}(k')\rangle\\
=& \langle g,\alpha^{-1}(h_{112})\rangle\langle f,\alpha^{-4}(k_{2})\alpha^{-3}(k')\rangle\langle g',h_{2}[S^{-1}(\alpha^{-5}(k_{1})\alpha^{-3}(h_{12}))\alpha^{-3}(h_{111})]\rangle\\
=& \langle g,h_{12}\rangle\langle f,\alpha^{-4}(k_{2})\alpha^{-3}(k')\rangle\\
&\langle g',[\alpha^{-2}(h_{22})S^{-1}(\alpha^{-2}(h_{21}))][S^{-1}(\alpha^{-4}(k_{1}))\alpha^{-2}(h_{11})]\rangle\\
=& \langle g,\alpha(h_{2})\rangle\langle f,\alpha^{-4}(k_{2})\alpha^{-3}(k')\rangle\langle g',[S^{-1}(\alpha^{-3}(k_{1}))h_{1}]\rangle.
\end{aligned}$$
Hence $\Delta^{op}(h\o f)R=R\Delta(h\o f)$.
Similarly we have $(\Delta\otimes\alpha)R=R^{13}R^{23},$
$(\alpha\otimes\Delta)R=R^{13}R^{12}.$

This finishes the proof. \hfill $\blacksquare$
\\

{\bf Example 2.20.}
In the Example 2.17, the Drinfeld double of $(kG,\phi)$ is given where $G$ is a finite group and $\phi$ is a group automorphism of $G$. Then by the above  proposition, the quasitriangular structure of $D(kG)$ is
$$R=\sum_{g\in G}1_{G}\bowtie e_{\phi(g)}\o g^{-1}\bowtie 1_{k^{G}}.$$
Especially when $G$ is a cycle group of order $n$, the quasitriangular structure of $D(kG)$ is
$$R=\sum_{0\leq i<n}1_{G}\bowtie e_{n-i}\o g^{-i}\bowtie 1_{k^{G}}.$$

\section*{3. DUAL PAIRS OF  HOM-HOPF-ALGEBRAS}
\def\theequation{3. \arabic{equation}}
\setcounter{equation} {0} \hskip\parindent

In this section, we will consider dual pairs of  Hom-Hopf algebras and
then give the Drinfel'd double associated to a pairing of Hom-Hopf algebras in the setting
 of dual pairs.
 \\

{\bf Definition 3.1.}
Let $(A,\alpha_{A})$ and $(B,\alpha_{B})$ be two Hom-Hopf algebras and $(-,-):A\o B\lr k$
 a non-degenerate bilinear form.  Then $(A,B)$ is called a {\sl dual pair} if the following conditions hold:
 \begin{eqnarray}
&&  (a,1)=\varepsilon(a),\quad \ (1,b)=\varepsilon(b),\\
&& (\alpha_{A}(a),\alpha_{B}(b))=(a,b), \\
&& (aa',b)=(\alpha^{2}_{A}(a),b_{1})(\alpha^{2}_{A}(a'),b_{2}),\\
&& (a,bb')=(a_{1},\alpha^{2}_{B}(b))(a_{2},\alpha^{2}_{B}(b')),\\
&& (S_{A}(a),b)=(a,S^{-1}_{B}(b)),
\end{eqnarray}
for any $a, a'\in A$ and $b, b'\in B$.
\\

{\bf Example 3.2.} Let $(A,B,\theta)$ be a dual pair of Hopf algebras defined in \cite{W}, $\alpha_{A}$ and $\alpha_{B}$ be a Hopf algebra isomorphism of $A$ and $B$ respectively. Assume that $(\alpha_{A},\alpha_{B},\theta)$ is a compatible pairing in the sense that
$$\theta(\alpha_{A}(a),\alpha_{B}(b))=\theta(a,b),
$$
for any $a\in A$ and $b\in B$. Define a bilinear form on $A_{\alpha_{A}}\o B_{\alpha_{B}}$ by
$$(\alpha_{A}(a),\alpha_{B}(b))=\theta(a,b),
$$
then $A_{\alpha_{A}}$ and $B_{\alpha_{B}}$ is a dual pair of Hom-Hopf algebra.
\\

Recall from \cite{M2} that a Hom-twisting map between the Hom-associative algebras $(A,\mu_{A},\alpha_{A})$ and $(B,\mu_{B},\alpha_{B})$ is a linear map $T: B\o A\lr A\o B$ such that the following conditions are satisfied:
\begin{eqnarray*}
&&  (\alpha_{A}\o\alpha_{B})\circ R=R\circ(\alpha_{B}\o\alpha_{A}),\\
&&  T\circ(\mu_{B}\o\alpha_{A})=(\alpha_{A}\o\mu_{B})\circ(T\o id)\circ(id\o T), \\
&&  T\circ(\alpha_{B}\o\mu_{A})=(\mu_{A}\o\alpha_{B})\circ(id\o T)\circ(T\o id).
\end{eqnarray*}

Then we will construct a Hom-twisting map. In fact, we have two linear maps $R_{1},R_{2}:A\o B\lr A\o B$ given by
$$R_{1}(a\o b)=(\alpha_{A}(a_{2}),b_{1})\alpha^{-1}_{A}(a_{1})\o\alpha^{-1}_{B}(b_{2}),$$
$$R_{2}(a\o b)=(\alpha_{A}(a_{1}),b_{2})\alpha^{-1}_{A}(a_{2})\o\alpha^{-1}_{B}(b_{1}).$$

Easy to check that both are invertible with the inverses
$$R^{-1}_{1}(a\o b)=(S^{-1}\alpha_{A}(a_{2}),b_{1})\alpha^{-1}_{A}(a_{1})\o\alpha^{-1}_{B}(b_{2}),$$
$$R^{-1}_{2}(a\o b)=(S^{-1}\alpha_{A}(a_{1}),b_{2})\alpha^{-1}_{A}(a_{2})\o\alpha^{-1}_{B}(b_{1}).$$
\\

Define the linear map $T=R_{1}\circ R^{-1}_{2}\circ\tau: B\o A\lr A\o B$, where $\tau$ is the flip map.
\\

{\bf Proposition 3.3.} The above $T$ is a Hom-twisting map.

{\bf Proof.}
Obviously $(\alpha_{A}\o\alpha_{B})\circ T=T\circ(\alpha_{B}\o\alpha_{A})$. Then for any $a,a'\in A$ and $b\in B$,
$$\begin{aligned}
&(\mu_{A}\o\alpha_{B})(id\o T)(T(b\o a)\o a')\\
=&\alpha^{-2}_{A}(a_{21}a'_{21})\o\alpha^{-2}_{B}(b_{221})( S^{-1}\alpha^{-1}_{A}(a_{1}a'_{1}),\alpha^{-3}_{B}(b_{222}))( a'_{22},\alpha^{-1}_{B}(b_{211}))( a_{22},b_{1})\\
=&\alpha^{-2}_{A}(a_{21}a'_{21})\o\alpha^{-1}_{B}(b_{21})( S^{-1}\alpha^{-1}_{A}(a_{1}a'_{1}),\alpha^{-2}_{B}(b_{22}))( a_{22}a'_{22},\alpha^{-1}_{B}(b_{1}))\\
=&\alpha^{-2}_{A}(a_{21}a'_{21})\o\alpha^{-1}_{B}(b_{12})(S^{-1}(a_{1}a'_{1}),b_{2})(a_{22}a'_{22},b_{11})\\
=&T(\alpha_{B}(b)\o aa').
\end{aligned}$$

Similarly we have $T\circ(\mu_{B}\o\alpha_{A})=(\alpha_{A}\o\mu_{B})\circ(T\o id)(id\o T).$ Therefore $T$ is a Hom-twisting map.

This finishes the proof. \hfill $\blacksquare$
\\

Hence we have the Hom-associative algebra $A\bowtie B$ with the multiplication
$$(a\bowtie b)(a'\bowtie b')=(S^{-1}\alpha_{A}(a'_{1}),b_{2})(a'_{22},\alpha^{-1}_{B}(b_{11}))a\alpha^{-2}_{A}(a'_{21})\bowtie\alpha^{-2}_{B}(b_{12})b'.$$
\\

{\bf Theorem 3.4.}
If we define $\Delta:A\bowtie B\lr A\bowtie B\o A\bowtie B$ by
$$\Delta(a\bowtie b)=a_{1}\bowtie b_{2}\o a_{2}\bowtie b_{1},$$
the antipode by
$$S=T\circ\tau\circ(S_{A}\o S^{-1}_{B}).$$
Then $A\bowtie B$ is a Hom-Hopf algebra, which is called the Drinfeld double of $A$ and $B$.
\\

{\bf Remark 3.5.} (1) Note that in the our construction, we don't need the Hom-Hopf algebra to be finite dimensional.

 (2) In fact for any finite dimensional Hom-Hopf algebra $(H,\alpha)$, $(H^{op},H^{*})$ is a dual pair if for any $h\in H$ and $f\in H^{*}$,
$(a,f)=\langle f,a\rangle.$ One can see that the Drinfeld double in this section is really a generalization of that in Corollary 2.16.
\\

{\bf Proposition 3.6.} By the multiplication, we can see that $A$ and $B$ could be imbedded into $A\bowtie B$ by $A\hookrightarrow A\bowtie B,\ a\mapsto a\bowtie1_{B}$ and $B\hookrightarrow A\bowtie B,\ b\mapsto 1_{A}\bowtie b$, respectively. Moreover $a\bowtie b=(\alpha^{-1}_{A}\o\alpha^{-1}_{B})((a\bowtie 1_{B})(1_{A}\bowtie b))$.
\\

\section*{4. THE DRINFEL'D DOUBLE VERSUS THE HEISENBERG DOUBLE  FOR HOM-HOPF ALGEBRAS}
\def\theequation{1. \arabic{equation}}
\setcounter{equation} {0} \hskip\parindent

In this section, the relation between the Drinfel'd double and Heisenberg double is established.
\\

{\bf Definition 4.1.} Let $(H,\alpha)$ be a Hom-Hopf algebra. The linear map $\sigma:H\o H\lr k$ is called a {\sl left Hom-2-cocycle}
 if the following conditions are satisfied

$$\sigma\circ(\alpha\o\alpha)=\sigma,\eqno(4.1)
$$

$$\sigma(l_{1},k_{1})\sigma(\alpha^{2}(h),l_{2}k_{2})=\sigma(h_{1},l_{1})\sigma(h_{2}l_{2},\alpha^{2}(k)),\eqno(4.2)
$$
for any $h,l,k\in H.$

Furthermore,  $\sigma$ is {\sl normal} if $\sigma(1,h)=\sigma(h,1)=\varepsilon(h).$
\\

Similarly if the condition (4.2) is replaced by
$$\sigma(\alpha^{2}(h),l_{1}k_{1})\sigma(l_{2},k_{2})=\sigma(h_{1}l_{1},\alpha^{2}(k))\sigma(h_{2},l_{2}),$$
then $\sigma$ is a {\sl  right Hom-2-cocycle}.
\\

{\bf Example 4.2.} (1) Recall from Corollary 2.15 that the Drinfeld double $D(H)$ of $(H,\alpha_{H})$  is the space $H^{op}\o H^{*}$ with the multiplication
$$(h\o f)(k\o g)= \alpha^{-2}_{H}(k_{21})h\o[\alpha^{-3}_{H}(k_{22})\rightharpoonup((\alpha^{*}_{H})^{2}(f)\leftharpoonup S\alpha^{-3}_{H}(k_{1}))]g,$$
for any $h,k\in H$ and $f,g\in H^{*}$.

 Define $\sigma:D(H)\o D(H)\lr k$ by
$$
\sigma(h\o f,k\o g)=\varepsilon(h)g(1)\langle f,\alpha(k)\rangle \eqno(4.3)
$$
Then it is not hard to verify that $\sigma$ is a left Hom-2-cocycle on $D(H)$.

(2) In \cite{MP}  Makhlouf and  Panaite  introduced another form of the Drinfeld double
 $\widetilde{D(A)}$, which is the space $(A^{op})^{*}\o A$ with the multiplication:
$$
(f\o a)(g\o b)=f[(\alpha^{-3}_{A}(a_{1})\rightharpoonup(\alpha^{2}_{A}(g)))
 \leftharpoonup S^{-1}\alpha^{-3}_{A}(a_{22})]\o\alpha^{-2}_{A}(a_{21})b,
$$
for any $a,b\in A$ and $f,g\in A^{*}$.

Define $\eta:\widetilde{D(A)}\o \widetilde{D(A)}\lr k$ by
$$
\eta(f\o a,g\o b)=\varepsilon(b)f(1)\langle g,\alpha(a)\rangle \eqno(4.4)
$$
Then $\eta$ is a right Hom-2-cocycle on $\widetilde{D(A)}$.
\\

{\bf Proposition 4.3.} Let $(H,\alpha)$ be a Hom-Hopf algebra.

(1) If $\sigma$ is a left Hom-2-cocycle, for any $h,k\in H$,
 define multiplication on $H$ as follows
$$h\cdot_{\sigma}k=\sigma(h_{1},k_{1})\alpha^{-1}(h_{2}k_{2}),$$
then $(H,\cdot_{\sigma},\alpha)$ is a Hom-associative algebra, called the left twist of $H$
 and denoted by ${}_{\si}H$.

(2) If $\sigma$ is a right Hom-2-cocycle, define multiplication on $H$ as follows
$$h\ _{\sigma}\cdot k=\alpha^{-1}(h_{1}k_{1})\sigma(h_{2},k_{2}),$$
then $(H,\cdot_{\sigma},\alpha)$ is also a Hom-associative algebra, called the right twist of $H$
 and denoted by $H_{\si}$.

{\bf Proof.} (1)  First of all, for any
$h,k,l\in H$, we have
$$1\cdot_{\sigma}h=h\cdot_{\sigma}1=\alpha(h),$$
and
$$\alpha(h\cdot_{\sigma}k)=\sigma(h_{1},k_{1})\alpha(\alpha^{-1}(h_{2}k_{2}))\stackrel{(4.1)}{ =}\alpha(h)\cdot_{\sigma}\alpha(k).$$
Then
\begin{eqnarray*}
\alpha(h)\cdot_{\sigma}(k\cdot_{\sigma}l)&=&\alpha(h)\cdot_{\sigma}(\sigma(k_{1},l_{1})\alpha^{-1}(k_{2}l_{2}))\\
                                         &=&\sigma(k_{11},l_{11})\sigma(\alpha^{2}(h_{1}), k_{12}l_{12}) h_{2}\alpha^{-1}(k_{2}l_{2})\\
                                         &\stackrel{(4.2)}{=}&\sigma(h_{11},k_{11})\sigma(h_{12}k_{12}, \alpha^{2}(l_{1}))\alpha^{-1}(h_{2}k_{2}) l_{2}\\
                                         &=&(h\cdot_{\sigma}k)\cdot_{\sigma}\alpha(l).
\end{eqnarray*}

(2) Similar to the proof of the part (1).

This finishes the proof. \hfill $\blacksquare$
\\

{\bf Definition 4.4.} The Heisenberg double  $H(A)$ of a Hom-Hopf algebra $(A,\alpha_{A})$ is the smash product $A\# A^{*}$
 with respect to the left regular action of $A^{*}$ on $A$, that is, it possesses the product as follows:
$$
(a\#f)(b\#g)=a(f_{1}\circ\alpha^{2}_{A}\rightharpoonup \alpha^{-1}_{A}(b))\#(f_{2}\circ\alpha_{A})g,
$$
for any $a,b\in A$ and $f,g\in A^{*}$.
\\

We now have the following main result of this section.
\\

{\bf Theorem 4.5.} Let  $(A,\alpha_{A})$ be a Hom-Hopf algebra.

  (1) The Heisenberg double $H(A^{op})$ of $A^{op}$ is the left twist of the Drinfeld double $D(A)$
   with the Hom-2-cocycle $\si$ given by Eq.(4.3).

  (2) The Heisenberg double $H(A^{*})$ of $A^{*}$ is the right twist of the Drinfeld double $\widetilde{D(A)}$
  with the right Hom-2-cocycle $\eta $ given by Eq.(4.4).

{\bf Proof.} (1) We need to check that ${}_{\si}D(A)$ and $H(A^{op})$ share the same multiplication. Indeed for any $f,g\in A^{*}$ and $a,b\in A$,
\begin{eqnarray*}
(a\o f)\cdot_{\sigma}(b\o g)
&=&\sigma(a_{1}\o f_{1},b_{1}\o g_{1})(\alpha^{-1}\o\alpha^{*})((a_{2}\o f_{2})(b_{2}\o g_{2}))\\
&=&\varepsilon(a_{1})\alpha^{-3}(b_{221})\alpha^{-1}(a_{2})\langle f_{1},\alpha(b_{1})\rangle\\
&&\quad \quad \o[\alpha^{-4}(b_{222})\rightharpoonup((\alpha^{3})^{*}(f_{2})\leftharpoonup S\alpha^{-4}(b_{21}))]g\\
&=&\alpha^{-1}(b_{1})a\o\langle f,\alpha^{-1}(?b_{2})\rangle g\\
&=&(a\#f)(b\#g)
\end{eqnarray*}

(2)  We just compute the multiplication in  $\widetilde{D(A)}$. For any $f,g\in A^{*}$ and $a,b\in A$,
$$\begin{aligned}
(f\o a)\cdot_{\eta}(g\o b)&=(\alpha^{*}\o\alpha^{-1})((f_{1}\o a_{1})(g_{1}\o b_{1}))\eta(f_{2}\o a_{2},g_{2}\o b_{2})\\
                            &=f[(\alpha^{-4}(a_{11})\rightharpoonup(\alpha^{3})^{*}(g_{2}))\leftharpoonup S^{-1}\alpha^{-4}(a_{122})]\\
                            &\quad \quad \o\alpha^{-3}(a_{121})b\langle f_{1},\alpha(a_{2})\rangle\\
                            &=f\langle g,\alpha^{-1}(?a_{1})\rangle\o \alpha^{-1}(a_{2})b,
\end{aligned}$$
which is exactly the multiplication in $H(A^{*})$.

This completes the proof. \hfill $\blacksquare$
\\

{\bf Example 4.6.} Let $(kG,\phi)$ be the Hom-Hopf algebra in Example 2.17 (2). By Example 4.2 (1), we have a left Hom-2-cocycle $\sigma$ on $D(kG)\o D(kG)$ given by
$$\sigma(g^{i}\o e_{j},g^{m}\o e_{k})=\delta_{k,0}\delta_{j,n-m},
$$
for any $0\leq i,j,m,k<n$. Then
\begin{eqnarray*}
&&(g^{i}\o e_{j})\cdot_{\sigma}(g^{m}\o e_{k})\\
&=&\sigma(g^{n-i}\o (e_{j})_{1},g^{n-m}\o (e_{k})_{1})(\phi^{-1}\o\phi^{*})[(g^{n-i}\o (e_{j})_{2})(g^{n-m}\o (e_{k})_{2})]\\
&=&(g^{i}\o e_{j+m})(g^{m}\o e_{k})=g^{n-(i+m)}\o \delta_{j+m,k}e_{n-k}\\
&=&(g^{i}\#e_{j})(g^{m}\# e_{k}),
\end{eqnarray*}
just as shown in Theorem 4.5.
\\

The following property is straightforward.
\\

 {\bf Proposition 4.7.}  Assume that $\sigma$ is a left Hom-2-cocycle on the Hom-Hopf algebra $(H,\alpha)$, then the comultiplication $\Delta$ of $H$ makes $_{\sigma}H$ into a right $H$-comodule Hom-algebra. Similar result holds for the algebra $H_{\sigma}$ if $\sigma$ is a right Hom-2-cocycle.
\\

We now can apply Theorem 4.5 and Proposition 4.7 to obtain the following corollary:
\\

{\bf Corollary 4.8.} (1) The comultiplication of $D(A)$, considered as a map from $H(A^{op})$ to $H(A^{op})\o D(A)$,  makes $H(A^{op})$ into a right $D(A)$-comodule Hom-algebra.

(2) The comultiplication of $\widetilde{D(A)}$, considered as a map from $H(A^{*})$ to $\widetilde{D(A)}\o H(A^{*})$, makes $H(A^{*})$ into a left $\widetilde{D(A)}$-comodule Hom-algebra.

\section*{ACKNOWLEDGEMENTS}

This work was supported by the NSF of China (No. 11371088) and the NSF of Jiangsu Province (No. BK2012736).

\end{document}